\documentclass[preprint,sort&compress,3p]{elsarticle}
\usepackage[utf8]{inputenc}
\usepackage{bbm}
\usepackage{amsmath}
\usepackage{amssymb}
\usepackage{amsthm}
\usepackage{graphicx}

\def\dualsymb{\circ}
\def\transposesymb{t}
\def\conjsymb{\ast}

\newcommand{\setC}{\mathbbm{C}}
\newcommand{\setN}{\mathbbm{N}}
\newcommand{\setR}{\mathbbm{R}}

\newcommand{\bounded}[1]{\mathcal{B}\left(#1\right)}
\newcommand{\boundedplus}[1]{\bounded{#1}^+}
\newcommand{\bbounded}[2]{\mathcal{B}\left(#1,#2\right)}
\newcommand{\hilbertspaceone}{\mathcal{K}}
\newcommand{\hilbertspacetwo}{\mathcal{H}}
\newcommand{\kh}{\hilbertspaceone\otimes\hilbertspacetwo}

\newcommand{\bk}{\bounded{\hilbertspaceone}}
\newcommand{\bh}{\bounded{\hilbertspacetwo}}

\newcommand{\bkh}{\bounded{\kh}}
\newcommand{\bhplus}{\bh^+}
\newcommand{\bkplus}{\bk^+}
\newcommand{\mappingcone}{\mathcal{C}}
\newcommand{\innerpr}[2]{\left<#1,#2\right>}
\newcommand{\innerprtwo}[2]{\innerpr{#1}{#2}'}
\newcommand{\innerprthree}[2]{\innerpr{#1}{#2}''}
\newcommand{\seq}[3]{\left\{#1\right\}_{#2}^{#3}}
\newcommand{\Tr}{\mathop{\textnormal{Tr}}}
\newcommand{\Ad}{\mathop{\textnormal{Ad}}}
\newcommand{\rk}{\mathop{\textnormal{rk}}}
\newcommand{\transposed}[1]{#1^{\transposesymb}}
\newcommand{\conj}[1]{#1^{\conjsymb}}
\newcommand{\dual}[1]{#1^{\dualsymb}}
\newcommand{\ddual}[1]{#1^{\dualsymb\dualsymb}}
\newcommand{\convhull}{\mathop{\textnormal{convhull}}}
\newcommand{\Id}{\mathop{\textnormal{id}}}
\newcommand{\One}{\mathbbm{1}}
\newcommand{\proj}[1]{p_{#1}}
\newcommand{\diad}[1]{\left|#1\right>\left<#1\right|}
\newcommand{\Choimatr}[1]{C_{#1}}
\newcommand{\matrices}[2]{M_{#1}\left(#2\right)}
\newcommand{\transpose}{\transposesymb}
\newcommand{\HPmaps}{\mathcal{HP}}
\newcommand{\Pmaps}{\mathcal{P}}
\newcommand{\kPmaps}[1]{\Pmaps_{#1}}
\newcommand{\SPmaps}{\mathcal{SP}}
\newcommand{\kSPmaps}[1]{\SPmaps_{#1}}
\newcommand{\CPmaps}{\mathcal{CP}}

\newcommand{\Pmapsb}[1]{\Pmaps\left(#1\right)}
\newcommand{\kPmapsb}[2]{\kPmaps{#1}\left(#2\right)}

\newcommand{\kSPmapsb}[2]{\kSPmaps{#1}\left(#2\right)}
\newcommand{\CPmapsb}[1]{\CPmaps\left(#1\right)}
\newcommand{\HPmapsbb}[2]{\HPmaps\left(#1,#2\right)}
\newcommand{\Pmapsbb}[2]{\Pmaps\left(#1,#2\right)}
\newcommand{\kPmapsbb}[3]{\kPmaps{#1}\left(#2,#3\right)}
\newcommand{\SPmapsbb}[2]{\SPmaps\left(#1,#2\right)}
\newcommand{\kSPmapsbb}[3]{\kSPmaps{#1}\left(#2,#3\right)}
\newcommand{\CPmapsbb}[2]{\CPmaps\left(#1,#2\right)}

\newtheorem{theorem}{Theorem}
\newtheorem{lemma}{Lemma}

\newtheorem{example}{Example}

\newtheorem{proposition}{Proposition}

\newtheorem{property}{Property}

\newtheorem{question}{Question}

\begin{document}

\title{Cones with a Mapping Cone Symmetry in the Finite-Dimensional Case}
\author{Łukasz Skowronek}
\address{Instytut Fizyki im. Smoluchowskiego, Uniwersytet Jagielloński, Reymonta 4, 30-059 Kraków, Poland}
\ead{lukasz.skowronek@uj.edu.pl}

\begin{abstract}
In the finite-dimensional case, we present a new approach to the theory of cones with a mapping cone symmetry, first introduced by St\o rmer. Our method is based on a definition of an inner product in the space of linear maps between two algebras of operators and the fact that the Jamio\l kowski-Choi isomorphism is an isometry. We consider a slightly modified class of cones, although not substantially different from the original mapping cones by St{\o}rmer. Using the new approach, several known results are proved faster and often in more generality than before. For example, the dual of a mapping cone turns out to be a mapping cone as well, without any additional assumptions. The main result of the paper is a characterization of cones with a mapping cone symmetry, saying that a given map is an element of such cone if and only if the composition of the map with the conjugate of an arbitrary element in the dual cone is completely positive. A similar result was known in the case where the map goes from an algebra of operators into itself and the cone is a symmetric mapping cone. Our result is proved without the additional assumptions of symmetry and equality between the domain and the target space. We show how it gives a number of older results as a corollary, including an exemplary application. 
\end{abstract}


\maketitle

\section{Introduction}
Mapping cones were introduced by St\o rmer in \cite{ref.St86} as a way to better understand the structure of positive maps. They are an abstract notion mimicking a well-known property of the cones of positive and more generally, $k$-positive maps. Namely, for any $k$-positive map $\Phi$ and a pair of completely positive maps $\Upsilon$ and $\Omega$, the map $\Upsilon\circ\Phi\circ\Omega$ is $k$-positive again.  St\o rmer called a closed cone $\mappingcone$, different from $0$, in the space of positive linear maps from an algebra $\bk$ of bounded operators on a Hilbert space $\hilbertspaceone$ into itself 
a \textit{mapping cone} \cite{ref.St86} if and only if for all $\Phi\in\mappingcone$ and $a,b\in\bk$, the map $x\mapsto a\Phi\left(bx\conj{b}\right)\conj{a}$ is an element of $\mappingcone$ again. Equivalently in the case of finite-dimensional $\hilbertspaceone$, $\Upsilon\circ\Phi\circ\Omega$ is an element of $\mappingcone$ for arbitrary completely positive maps $\Upsilon,\Omega$  of $\bk$. In contrast to the original paper \cite{ref.St86}, in the following we do not only consider cones of maps from $\bk$ into itself, but also into $\bh$ for another Hilbert space $\hilbertspacetwo$. Most of the time, we use the convexity assumption, which was absent in \cite{ref.St86}, but is actually very much in line of St\o rmer's later work on mapping cones. Even though the mentioned differences are not substantial, we shall stick to the term ``cones with a mapping cone symmetry'' (\textit{mcs-cones} for short) in order to give sufficient credit to \cite{ref.St86}.

In a number of recent papers \cite{ref.St08,ref.St09dual,ref.St09sepposII, ref.SSZ09}, St\o rmer and coauthors proved various characterization theorems for mapping cones. Some of them are of special interest because they relate purely geometrical properties to properties of algebraic nature. In particular, they reveal an intrinsic link between the condition that a product of two maps is completely positive and the fact that the two maps belong to a pair of dual mapping cones (cf. \cite{ref.St09dual} and \cite{ref.SSZ09}). In the present paper, we aim at a similar, and very general characterization for mcs-cones. We also use a new approach that allows us to prove results more quickly and to directly exploit the mapping cone symmetry. Our methods work well in the finite-dimensional setting, whereas their applicability to the infinite-dimensional case is not obvious.

\section{Basic notions}
Let $\hilbertspaceone$ and $\hilbertspacetwo$ be two Hilbert spaces. We denote with $\innerpr{.}{.}$ the inner product in $\hilbertspaceone$ or $\hilbertspacetwo$. In the following, we shall assume that $\hilbertspaceone$ and $\hilbertspacetwo$ are finite-dimensional and thus equivalent to $\setC^m$ and $\setC^n$ for some $m,n\in\setN$, $\dim\hilbertspaceone=m$, $\dim\hilbertspacetwo=n$. We also fix orthonormal bases $\seq{f_j}{j=1}{m}$ and $\seq{e_i}{i=1}{n}$ of $\hilbertspaceone$ and $\hilbertspacetwo$, resp. Thus we have a very specific setting for our discussion, but we shall keep the abstract notation of Hilbert spaces, hoping to bring the attention of the reader to possible generalizations to the infinite-dimensional case. Let us denote with $\bk$ and $\bh$ the spaces of bounded operators on $\hilbertspaceone$ and $\hilbertspacetwo$ resp. and choose their canonical bases $\seq{f_{kl}}{k,l=1}{m}$, $\seq{e_{ij}}{i,j=1}{n}$. That is, $f_{kl}\left(e_j\right)=\delta_{lj}f_k$ and similarly for the $e_{ij}$.  Positive elements of $\bk$ are operators $A\in\bk$ such that $\innerpr{v}{A\left(v\right)}\geqslant 0\,\forall_{v\in\hilbertspacetwo}$. Similarly for elements of $\bh$. The sets of positive elements of $\bk$ and $\bh$ will be denoted by $\bkplus$ and $\bhplus$. In the finite-dimensional case, there exists a natural inner product in $\bk$, given by the formula
\begin{equation}\label{HSProd}
 \innerprtwo{A}{B}:=\Tr\left(A\conj{B}\right)
\end{equation}
for $A,B\in\bk$. An identical definition works for $A,B\in\bh$ and we do not distinguish notationally between the inner products in $\bh$ and $\bk$. Note that the bases $\seq{f_{kl}}{k,l=1}{m}$ and $\seq{e_{ij}}{i,j=1}{n}$ are orthonormal with respect to $\innerprtwo{.}{.}$.

In the following, we will be mostly dealing with linear maps from $\bk$ to $\bh$. Because of the finite-dimensionality assumption, they are all elements of $\bbounded{\bk}{\bh}$, the space of bounded operators from $\bk$ to $\bh$. Given a map $\Phi\in\bbounded{\bk}{\bh}$, we define its conjugate $\conj{\Phi}$ as a map from $\bh$ into $\bk$ satisfying $\innerprtwo{A}{\Phi\left(B\right)}=\innerprtwo{\conj{\Phi}\left(A\right)}{B}$ for all $A\in\bh$ and $B\in\bk$. In our setting, there also exists a natural inner product in $\bbounded{\bk}{\bh}$, given by the formula
\begin{equation}\label{HSProdtwo}
 \innerprthree{\Phi}{\Psi}:=\sum_{k,l=1}^{m}\innerprtwo{\Phi\left(f_{kl}\right)}{\Psi\left(f_{kl}\right)}.
\end{equation}
Note that the spaces $\bbounded{\bh}{\bk}$, $\bounded{\bk}$ and $\bounded{\bh}$ can be endowed with analogous inner products and we shall not notationally distinguish between them. The following proposition summarizes a few elementary facts about $\innerprthree{.}{.}$ that will be useful for our later discussion.

\begin{proposition}\label{propinnerpr}
 For all $\Phi,\Psi\in\bbounded{\bk}{\bh}$ and $\alpha\in\bounded{\bh}$, $\beta\in\bounded{\bk}$, one has the following equalities
\begin{enumerate}
\item $\innerprthree{\Phi\circ\beta}{\Psi}=\innerprthree{\beta}{\conj{\Phi}\circ\Psi}=\innerprthree{\conj{\Psi}\circ\Phi}{\conj{\beta}}$,
\item $\innerprthree{\alpha\circ\Phi}{\Psi}=\innerprthree{\alpha}{\Psi\circ\conj{\Phi}}=\innerprthree{\Phi\circ\conj{\Psi}}{\conj{\alpha}}$,
\item $\innerprthree{\alpha\circ\Phi\circ\beta}{\Psi}=\innerprthree{\Phi}{\conj{\alpha}\circ\Psi\circ\conj{\beta}}$.
\end{enumerate}

\begin{proof} The first equality in point one follows  directly from $\innerprtwo{\Phi\circ\beta\left(f_{kl}\right)}{\Psi\left(f_{kl}\right)}=\innerprtwo{\beta\left(f_{kl}\right)}{\conj{\Phi}\circ\Psi\left(f_{kl}\right)}$ and the definition of $\innerprthree{.}{.}$, eq. \eqref{HSProdtwo}. To prove the other equalities, we can use a simple lemma.
\begin{lemma}\label{lemmapoema}
 For any finite-dimensional Hilbert spaces $\hilbertspaceone$, $\hilbertspacetwo$ and maps $\Phi,\Psi\in\bbounded{\bk}{\bh}$, we have
\begin{equation}
\innerprthree{\Phi}{\Psi}=\innerprthree{\conj{\Psi}}{\conj{\Phi}} 
\end{equation}
\begin{proof} Starting from the definition of $\innerprthree{.}{.}$, we get
\begin{multline}\label{proofLemma1}
 \innerprthree{\Phi}{\Psi}=\sum_{k,l=1}^m\innerprtwo{\Phi\left(f_{kl}\right)}{\Psi\left(f_{kl}\right)}=\sum_{i,j=1}^n\sum_{m,n=1}^n\sum_{k,l=1}^m\Phi_{ij,kl}\overline{\Psi_{mn,kl}}\innerprtwo{e_{ij}}{e_{mn}}=\sum_{i,j=1}^n\sum_{k,l=1}^m\Phi_{ij,kl}\overline{\Psi_{ij,kl}}=\\
=\sum_{i,j=1}^n\sum_{k,l=1}^m\sum_{r,s=1}^m\Phi_{ij,rs}\overline{\Psi_{ij,kl}}\innerprtwo{f_{rs}}{f_{kl}}=\sum_{i,j=1}^n\sum_{k,l=1}^m\sum_{r,s=1}^m\innerprtwo{\overline{\Psi_{ij,rs}}f_{r,s}}{\overline{\Phi_{ij,kl}}f_{kl}}=\sum_{i,j=1}^n\innerprtwo{\conj{\Psi}\left(e_{ij}\right)}{\conj{\Phi}\left(e_{ij}\right)},
\end{multline}
where the last equality follows because $\conj{\Phi}\left(e_{ij}\right)=\sum_{k,l=1}^m\overline{\Phi_{ij,kl}}f_{kl}$ as a consequence of $\innerprtwo{f_{kl}}{\conj{\Phi}\left(e_{ij}\right)}=\innerprtwo{\Phi\left(f_{kl}\right)}{e_{ij}}=\sum_{r,s=1}^m\overline{\Phi_{rs,kl}}\innerprtwo{e_{rs}}{e_{ij}}=\overline{\Phi_{ij,kl}}$. Similarly, $\conj{\Psi}\left(e_{ij}\right)=\sum_{r,s=1}^m\overline{\Phi_{ij,rs}}f_{rs}$ holds.
The final expression in \eqref{proofLemma1} clearly equals $\innerprthree{\conj{\Psi}}{\conj{\Phi}}$.
\end{proof} 
\end{lemma}
Note that the assertion of Lemma \ref{lemmapoema} holds for any choice of $\hilbertspaceone$ and  $\hilbertspacetwo$, and thus also when the two finite-dimensional Hilbert spaces are different from the $\hilbertspaceone$ and  $\hilbertspacetwo$ referred to in the statement of the proposition. Using the lemma, we get $\innerprthree{\beta}{\conj{\Phi}\circ\Psi}=\innerprthree{\conj{\Psi}\circ\Phi}{\conj{\beta}}$, which proves the second equality in point one. Furthermore,
\begin{equation}\label{eqpointtwo}
 \innerprthree{\alpha\circ\Phi}{\Psi}=\innerprthree{\conj{\Psi}}{\conj{\Phi}\circ\conj{\alpha}}=\overline{\innerprthree{\conj{\Phi}\circ\conj{\alpha}}{\conj{\Psi}}}=\overline{\innerprthree{\conj{\alpha}}{\Phi\circ\conj{\Psi}}}=\innerprthree{\Phi\circ\conj{\Psi}}{\conj{\alpha}}=\innerprthree{\alpha}{\Psi\circ\conj{\Phi}},
\end{equation}
where we successively used Lemma \ref{lemmapoema}, the conjugate symmetry of $\innerprthree{.}{.}$, the first equation in point one, the conjugate symmetry again, and finally Lemma \ref{lemmapoema} for the second time. Obviously, the first, the fifth and the sixth term in equation \eqref{eqpointtwo} are the same as in point two of the proposition. Hence the only remaining thing to prove is point three. We have
\begin{equation}\label{eqpointthree} 
 \innerprthree{\alpha\circ\Phi\circ\beta}{\Psi}=\innerprthree{\alpha}{\Psi\circ\conj{\beta}\circ\conj{\Phi}}=\innerprthree{\beta\circ\conj{\Psi}\circ\alpha}{\conj{\Phi}}=\innerprthree{\Phi}{\conj{\alpha}\circ\Psi\circ\conj{\beta}},
\end{equation}
where we used the properties $\innerprthree{\alpha\circ\Phi}{\Psi}=\innerprthree{\alpha}{\Psi\circ\conj{\Phi}}$ with $\Phi\rightarrow\Phi\circ\beta$, $\innerprthree{\beta}{\conj{\Phi}\circ\Psi}=\innerprthree{\Phi\circ\beta}{\Psi}$ with $\beta\rightarrow\alpha$, $\Phi\rightarrow\beta\circ\conj{\Psi}$ and $\Psi\rightarrow\conj{\Phi}$, and finally Lemma \ref{lemmapoema}.
 \end{proof}
\end{proposition}

Consider the tensor product $\kh$. This space has a natural inner product, inherited from $\hilbertspaceone$ and $\hilbertspacetwo$, and an orthonormal basis $\seq{f_{kl}\otimes e_{ij}}{i,j=1;k,l=1}{n;m}$. Similarly to $\bk$ and $\bh$, the space $\bkh$ of bounded operators on $\kh$ is endowed with a natural Hilbert-Schmidt product, defined by formula \eqref{HSProd} with $A,B\in\bkh$. We shall again denote the inner product with $\innerprtwo{.}{.}$ to avoid excess notation. There exists a one-to-one correspondence between linear maps $\Phi$ of $\bk$ into $\bh$ and elements of $\bkh$, given by
\begin{equation}\label{Jamisodef}
\Phi\mapsto\Choimatr{\Phi}:=\sum_{k,l=1}^m f_{kl}\otimes\Phi\left(f_{kl}\right).
\end{equation}
The symbol $\Choimatr{\Phi}$ denotes the {\it Choi matrix} of $\Phi$ \cite{ref.Choi75} and the mapping $J:\Phi\mapsto\Choimatr{\Phi}$ is sometimes called the {\it Jamio\l kowski-Choi isomorphism} \cite{ref.J72}. In fact, $J$ is not only an isomorphism, but also an \emph{isometry} between $\bbounded{\bk}{\bh}$ and $\bkh$ in the sense of Hilbert-Schmidt type inner products. One has the following
\begin{property}\label{propertythree}The Jamio\l kowski-Choi isomorphism is an isometry. One has
\begin{equation}
 \innerprthree{\Phi}{\Psi}=\innerprtwo{\Choimatr{\Phi}}{\Choimatr{\Psi}}
\end{equation}
for all $\Phi,\Psi\in\bbounded{\bk}{\bh}$ (with $\Choimatr{\Phi},\Choimatr{\Psi}\in\bkh$).
\begin{proof}By the definition of $\Choimatr{\Phi}$ and $\Choimatr{\Psi}$,
\begin{equation}\label{proofJamisoiso1}
 \innerprtwo{\Choimatr{\Phi}}{\Choimatr{\Psi}}=\innerprtwo{\sum_{k,l=1}^mf_{kl}\otimes\Phi\left(f_{kl}\right)}{\sum_{r,s=1}^mf_{rs}\otimes\Psi\left(f_{rs}\right)}=\ldots
\end{equation}
Since $\Tr\left(\left(A\otimes A'\right)\conj{\left(B\otimes B'\right)}\right)=\Tr\left(A\conj{B}\right)\Tr\left(A'\conj{B'}\right)$ for arbitrary $A,B\in\bk$ and $A',B'\in\bh$, by formula \eqref{HSProd} we have
\begin{equation}\label{proofJamisoiso2}
 \ldots=\sum_{k,l=1}^m\sum_{r,s=1}^m\innerprtwo{f_{kl}}{f_{rs}}\innerprtwo{\Phi\left(f_{kl}\right)}{\Psi\left(f_{rs}\right)}=\sum_{k,l=1}^m\innerprtwo{\Phi\left(f_{kl}\right)}{\Psi\left(f_{kl}\right)},
\end{equation}
where we used orthonormality of $\seq{f_{kl}}{k,l=1}{m}$.
The last expression equals $\innerprthree{\Phi}{\Psi}$ by definition \eqref{HSProdtwo}.
 \end{proof}
\end{property}

A linear map $\Phi$ from $\bk$ to $\bh$ is called {\it positive} if it preserves positivity of operators, which means $\Phi\left(\bkplus\right)\subset\bhplus$. Moreover, $\Phi$ is called {\it$k$-positive} if $\Phi\otimes\Id_{\matrices{k}{\setC}}$ is positive as a map from $\bk\otimes\matrices{k}{\setC}$ into $\bh\otimes\matrices{k}{\setC}$, where $\matrices{k}{\setC}$ denotes the space of $k\times k$ matrices with complex entries and $\Id$ refers to the identity map. A map $\Phi$ is called {\it completely positive} if it is $k$-positive for all $k\in\setN$. From the Choi's theorem on completely positive maps \cite{ref.Choi75} (cf. also Lemma \ref{lemmaCofAdV}) it follows that every such map has a representation $\Phi=\sum_i\Ad_{V_i}$ as a sum of {\it conjugation maps}, $\Ad_{V_i}:\rho\mapsto V_i\rho\conj{V_i}$ with $V_i\in\bbounded{\hilbertspaceone}{\hilbertspacetwo}$. Conversely, every map $\Phi$ of the form $\sum_i\Ad_{V_i}$ is completely positive. If all the $V_i$'s can be chosen of rank $\leqslant k$ for some $k\in\setN$, $\Phi$ is said to be {\it$k$-superpositive} \cite{ref.SSZ09}. One-superpositive maps are simply called {\it superpositive} \cite{ref.Ando04}. The sets of positive, $k$-positive, completely positive, $k$-superpositive and superpositive maps from $\bk$ to $\bh$ will be denoted with $\Pmapsbb{\bk}{\bh}$, $\kPmapsbb{k}{\bk}{\bh}$, $\CPmapsbb{\bk}{\bh}$, $\kSPmapsbb{k}{\bk}{\bh}$, $\SPmapsbb{\bk}{\bh}$ or $\Pmaps$, $\kPmaps{k}$, $\CPmaps$, $\kSPmaps{k}$, $\SPmaps$ for short. It is clear that all of them are closed convex cones contained in $\Pmapsbb{\bk}{\bh}$. They also share a more special property that the product $\Upsilon\circ\Phi\circ\Omega$ of $\Phi\in\mappingcone$, $\Upsilon\in\CPmapsb{\bh}$ and $\Omega\in\CPmapsb{\bk}$ is an element of $\mappingcone$ again, where $\mappingcone$ stands for one of the sets $\Pmaps$, $\kPmaps{k}$, $\CPmaps$, $\kSPmaps{k}$ and $\SPmaps$ (cf. e.g. \cite{ref.SSZ09}). Thus, following rather closely the original definition by St{\o}rmer \cite{ref.St86}, a {\it cone with a mapping cone symmetry}, or an {\it mcs-cone} for short, is defined as a closed convex cone $\mappingcone$ in $\Pmapsbb{\bk}{\bh}$, different from $0$, such that $\Upsilon\circ\Phi\circ\Omega\in\mappingcone$ for all $\Phi\in\mappingcone$, 
$\Upsilon\in\CPmapsb{\bh}$ and $\Omega\in\CPmapsb{\bk}$. In the following, the convexity assumption could sometimes be skept, and we do include appropriate comments. 

Note that the set of positive maps from $\bk$ into $\bh$ is contained in the real-linear subspace $\HPmapsbb{\bk}{\bh}\subset\bbounded{\bk}{\bh}$ ($\HPmaps$ for short) consisting of all Hermiticity-preserving maps, i.e. $\Phi$ such that $\Phi\left(\conj{X}\right)=\conj{\Phi\left(X\right)}$. Moreover, the image of $\HPmapsbb{\bk}{\bh}$ by $J:\Phi\mapsto\Choimatr{\Phi}$ equals the set of self-adjoint elements of $\bkh$ \cite{ref.Pillis}. Therefore $\innerprthree{.}{.}$ induces a \emph{symmetric}
 inner product on $\HPmapsbb{\bk}{\bh}$ (cf. Property \ref{propertythree}). By definition, all mapping cones are subsets of $\Pmaps$ and thus of $\HPmaps$. Since $\HPmaps$ is a finite-dimensional space over $\setR$ with a symmetric inner product $\innerprthree{.}{.}$, one can easily apply to it tools of convex analysis. In particular, given any cone $\mappingcone\subset\HPmaps$, one defines its {\it dual} $\dual{\mappingcone}$ as the cone of elements $\Psi\in\HPmaps$ such that $\innerprthree{\Psi}{\Phi}\geqslant 0$ for all $\Phi\in\mappingcone$,
\begin{equation}\label{dualconedef}
 \dual{\mappingcone}:=\left\{\Psi\in\HPmapsbb{\bk}{\bh}\,\vline\,\innerprthree{\Psi}{\Phi}\geqslant 0\,\forall_{\Phi\in\mappingcone}\right\}.
\end{equation}
 Obviously, $\dual{\mappingcone}$ is closed and convex. It has a clear geometrical interpretation as the convex cone spanned by the normals to the supporting hyperplanes for $\mappingcone$. The dual cone has a well-known counterpart in convex analysis \cite{ref.Rockafellar}, $\mappingcone^{\star}=-\dual{\mappingcone}$, which is called the {\it polar} of $\mappingcone$. We have the following
\begin{property}\label{propertyfour}
 Let $\mappingcone$ be a closed convex cone. Then
\begin{equation}\label{doubledual}
 \mappingcone=\ddual{\mappingcone}.
\end{equation}
\begin{proof}
 Formula \eqref{propertyfour} is equivalent to $\mappingcone^{\star\star}=\mappingcone$ for a closed convex cone $\mappingcone$. The latter equality is a known fact in convex analysis. A proof can be found e.g. in \cite{ref.Rockafellar} (Theorem 14.1). 
\end{proof}
\end{property}
It can be shown (cf. e.g. \cite{ref.SSZ09}) that a duality relation $\dual{\kPmaps{k}}=\kSPmaps{k}$ holds for all $k\in\setN$. The converse relation $\dual{\kSPmaps{k}}=\kPmaps{k}$ is also true, as a consequence of Property \ref{propertyfour}. In particular, for $k=1$ we get $\dual{\SPmaps}=\Pmaps$ and $\dual{\Pmaps}=\SPmaps$. Taking $k=\min\left\{m,n\right\}$, one obtains $\dual{\CPmaps}=\CPmaps$, which is in accordance with Choi's theorem on completely postive maps \cite{ref.Choi75} and with Property \ref{propertythree}.

In the following, we shall be interested in duality relations between mcs-cones. This is in general a well-posed problem, because the operation $\mappingcone\rightarrow\dual{\mappingcone}$ acts within the ``mcs'' class. We have
\begin{proposition}\label{dualisamappingcone}
Let $\mappingcone\subset\Pmapsbb{\bk}{\bh}$ be an arbitrary mcs-cone. Then $\dual{\mappingcone}$, defined as in \eqref{dualconedef}, is an mcs-cone as well.
\begin{proof}
 Let $\Psi$ be an element of $\dual{\mappingcone}$. First we prove that $\Upsilon\circ\Psi\circ\Omega\in\dual{\mappingcone}$ for all $\Upsilon\in\CPmapsb{\bh}$ and $\Omega\in\CPmapsb{\bk}$. We have $\conj{\Upsilon}\in\CPmapsb{\bh}$ and $\conj{\Omega}\in\CPmapsb{\bk}$ because the sets of completely positive maps are $\conjsymb$-invariant. Therefore $\conj{\Upsilon}\circ\Phi\circ\conj{\Omega}\in\mappingcone$ for an arbitrary element $\Phi$ of the cone $\mappingcone$. By the definition \eqref{dualconedef} of $\dual{\mappingcone}$, we have
$ \innerprthree{\Psi}{\conj{\Upsilon}\circ\Phi\circ\conj{\Omega}}\geqslant 0\,\forall_{\Phi\in\mappingcone}$. Using Proposition \ref{propinnerpr}, point three, we can rewrite this as
\begin{equation}\label{dualineq2}
 \innerprthree{\Upsilon\circ\Psi\circ\Omega}{\Phi}\geqslant 0\,\forall_{\Phi\in\mappingcone}.
\end{equation}
According to definition \eqref{dualconedef}, condition \eqref{dualineq2} means that $\Upsilon\circ\Psi\circ\Omega\in\dual{\mappingcone}$. This holds for arbitrary $\Upsilon\in\CPmapsb{\bh}$ and $\Omega\in\CPmapsb{\bk}$. The only thing which is left to prove is $\dual{\mappingcone}\subset\Pmapsbb{\bk}{\bh}$. The inclusion holds because every mcs-cone $\mappingcone$ contains all the conjugation maps $\Ad_V$ with $\rk V=1$. Consequently, $\dual{\mappingcone}\subset\dual{\convhull\left\{\Ad_V\,\vline\,\rk V=1\right\}}=\dual{\SPmaps}=\Pmaps$. To show that indeed $\left\{\Ad_V\vline\rk V=1\right\}\subset\mappingcone$ for any mcs-cone $\mappingcone$, take an arbitrary nonzero $\Phi\in\mappingcone$. There must exist normalized vectors $\upsilon\in\hilbertspaceone$ and $\omega\in\hilbertspacetwo$ such that $\innerprtwo{\proj{\omega}}{\Phi\left(\proj{\upsilon}\right)}\geqslant 0$, where $\proj{\upsilon}$ and $\proj{\omega}$ are orthogonal projections onto the one-dimensional subspaces spanned by $\upsilon$ and $\omega$. Denote $\chi:=\innerprtwo{\proj{\omega}}{\Phi\left(\proj{\upsilon}\right)}$. Consider a pair of maps, $U:\hilbertspaceone\ni a\mapsto\innerpr{a}{\upsilon'}\upsilon\in\hilbertspaceone$ and $W:\hilbertspacetwo\ni b\mapsto\innerpr{b}{\omega}\omega'\in\hilbertspacetwo$, where $\upsilon'$ and $\omega'$ are arbitrary normalized vectors in $\hilbertspaceone$ and $\hilbertspacetwo$. A map $\Phi'$, defined as $\lambda/\chi\left(\Ad_W\circ\,\Phi\circ\Ad_U\right)$ acts in the following
 way, $\Phi':\rho\mapsto\lambda\innerprtwo{\proj{\upsilon'}}{\rho}\proj{\omega'}$ or $\Phi'=\Ad_V$ with $V:\hilbertspaceone\ni c\mapsto\lambda\innerpr{\upsilon'}{c}\omega'$. Any rank one operator $V$ can be written in the latter form for some $\upsilon'$ and $\omega'$. But $\Phi'$ is an element of $\mappingcone$ because of the assumption that $\mappingcone$ is an mcs-cone. Thus indeed $\Ad_V\in\mappingcone$ for all $V\in\bbounded{\hilbertspaceone}{\hilbertspacetwo}$ such that $\rk V=1$. In the case of $\hilbertspaceone=\hilbertspacetwo$ and mapping cones $\mappingcone$ as in the original definition by St{\o}rmer,  the inclusion $\Ad_V\in\mappingcone$ follows from Lemma 2.4 in \cite{ref.St86}. It should be kept in mind that we never used convexity of $\mappingcone$ in the proof.
\end{proof}
\end{proposition}

Note that a version of Proposition \ref{dualisamappingcone} was proved in \cite{ref.St09mappingcones} using different methods, with the additional assumption of $\hilbertspacetwo=\hilbertspaceone$ and $\mappingcone$ being a symmetric mapping cone. It is instructive to see how that result of \cite{ref.St09mappingcones} follows using our method. First, note that a mapping cone $\mappingcone\subset\Pmapsb{\bk}$ is called {\it symmetric} \cite{ref.St09mappingcones} if $\mappingcone=\conj{\mappingcone}=\transposed{\mappingcone}$, with $\conj{\mappingcone}:=\left\{\conj{\Phi}|\Phi\in\mappingcone\right\}$ and $\transposed{\mappingcone}:=\left\{\transpose\circ\Phi\circ\transpose|\Phi\in\mappingcone\right\}$, where $\transpose$ stands for the transposition map, $t:\bk\ni f_{kl}\mapsto f_{lk}\in\bk$. We have the following
\begin{proposition}\label{dualsymmetricissymmetric}Consider the case $\hilbertspacetwo=\hilbertspaceone$.
 Let $\mappingcone$ be a symmetric mapping cone of maps from $\bk$ into itself, i.e. $\mappingcone\in\Pmapsb{\bk}$. The dual $\dual{\mappingcone}$ is also a symmetric mapping cone.
\begin{proof}
 By Proposition \ref{dualisamappingcone} and the redundancy of the convexity assumption, we know that $\dual{\mappingcone}$ is an mcs-cone of maps from $\bk$ into itself, and thus a mapping cone in the sense of \cite{ref.St86}. We only need to show that it is symmetric. Let $\Psi$ be an arbitrary element of $\dual{\mappingcone}$. By the symmetry $\transposed{\mappingcone}=\mappingcone$, we know that the condition $\innerprthree{\Psi}{\Phi}\geqslant 0\,\forall_{\Phi\in\mappingcone}$ is equivalent to $\innerprthree{\Psi}{\transpose\circ\Phi\circ\transpose}\geqslant 0\,\forall_{\Phi\in\mappingcone}$. By Proposition \ref{propinnerpr}, point three, this is the same as $\innerprthree{\transpose\circ\Psi\circ\transpose}{\Phi}\geqslant 0\,\forall_{\Phi\in\mappingcone}$, or $\transpose\circ\Psi\circ\transpose\in\dual{\mappingcone}$. Thus $\transposed{\left(\dual{\mappingcone}\right)}=\dual{\mappingcone}$. To show $\conj{\left(\dual{\mappingcone}\right)}=\dual{\mappingcone}$, one only needs to note that $\innerprthree{\Psi}{\Phi}=\innerprthree{\conj{\Phi}}{\conj{\Psi}}$ by Lemma \ref{lemmapoema}. Now, the property $\conj{\mappingcone}=\mappingcone$ can be used.
\end{proof}
\end{proposition}

\section{The main theorem}
Using the properties discussed in the previous section, we can almost immediately prove a surprising characterization theorem for mcs-cones, which was strongly suggested by earlier results on the subject \cite{ref.St09dual,ref.St09mappingcones,ref.SSZ09}. It holds without any additional assumptions about the cone, and is noteworthy as it links the condition that two maps $\Phi$, $\Psi$ lay in a pair of dual mcs-cones to the fact that the product $\conj{\Psi}\circ\Phi$ is a $\CPmaps$ map. Thus it reveals a connection between convex geometry and a fact which is more likely to be called algebraic than geometrical. Before we proceed with the proof, let us show a simple lemma, which is a version of \cite[Lemma 1$(i)$]{ref.SS10} for $\hilbertspaceone\neq\hilbertspacetwo$.
\begin{lemma}\label{lemmaCofAdV}Let $V:\hilbertspaceone\ni a\mapsto\sum_{i=1}^n\sum_{j=1}^mV_{ij}\innerpr{a}{f_j}e_i\in\hilbertspacetwo$ be an arbitrary operator in $\bbounded{\bk}{\bh}$ and consider the map $\Ad_V:\rho\mapsto V\rho\conj{V}$. Then
\begin{equation}\label{CAdV1}
\Choimatr{\Ad_V}=\diad{\upsilon}, 
\end{equation}
where $\upsilon=\sum_{i=1}^n\sum_{j=1}^mV_{ij}f_j\otimes e_i$ is a vector in $\hilbertspaceone\otimes\hilbertspacetwo$ and $\diad{\upsilon}:w\mapsto\innerpr{w}{\upsilon}\upsilon$ is proportional to an orthogonal projection onto the subspace spanned by $\upsilon$.
\begin{proof}Obviously, the map $\conj{V}$ acts in the following way, $\conj{V}:\hilbertspacetwo\ni b\mapsto\sum_{i=1}^n\sum_{j=1}^m\overline{V_{ij}}\innerpr{b}{e_i}f_j\in\hilbertspaceone$. Thus 
\begin{equation}\label{VfklV}
Vf_{kl}\conj{V}:\hilbertspacetwo\ni b\mapsto\sum_{i,r=1}^n\sum_{j,s=1}^mV_{rs}\innerpr{f_{kl}\left(f_j\right)}{f_{s}}\overline{V_{ij}}\innerpr{b}{e_i}e_{r}\in\hilbertspacetwo,
\end{equation}
where the last expression is easily verified to be equal to $\sum_{i,r=1}^nV_{rk}\overline{V_{il}}\innerpr{b}{e_i}e_{r}$. Thus we have $Vf_{kl}\conj{V}=\sum_{i,r=1}^nV_{rk}\overline{V_{il}}e_{ri}$ and by the definition \eqref{Jamisodef} of the Choi matrix,
\begin{equation}\label{CAdV2}
 \Choimatr{\Ad_V}=\sum_{k,l=1}^m\sum_{i,r=1}^nV_{rk}\overline{V_{il}}f_{kl}\otimes e_{ri}=\diad{\upsilon},
\end{equation}
with $\upsilon=\sum_{i=1}^n\sum_{j=1}^mV_{ij}f_j\otimes e_i$. A proof of the last equality in \eqref{CAdV2} is left as an elementary exercise for the reader. 
\end{proof}
\end{lemma}

We are ready to prove the following (cf. Theorem 1 in \cite{ref.St09dual}).

\begin{theorem}\label{MainTheorem}
 Let $\mappingcone\subset\Pmapsbb{\bk}{\bh}$ be an mcs-cone. The following conditions are equivalent,
\begin{enumerate}
 \item $\Phi\in\mappingcone$,
\item $\conj{\Psi}\circ\Phi\in\CPmapsb{\bk}$ for all $\Psi\in\dual{\mappingcone}$,
\item $\Phi\circ\conj{\Psi}\in\CPmapsb{\bh}$ for all $\Psi\in\dual{\mappingcone}$.
\end{enumerate}
\begin{proof}
 We first show $1\Leftrightarrow 2$. Let us start with $2\Rightarrow 1$. Since $\conj{\Psi}\circ\Phi\in\CPmaps\,\forall_{\Psi\in\dual{\mappingcone}}$, we can use the facts that $\dual{\CPmaps}=\CPmaps$ and $\Id\in\CPmaps$ to get 
\begin{equation}\label{conjPsiPhiid}
 \innerprthree{\conj{\Psi}\circ\Phi}{\Id}\geqslant 0\,\forall_{\Psi\in\dual{\mappingcone}}.
\end{equation}
By using point one of Proposition \ref{propinnerpr} with the identity map $\Id$ substituted for $\beta$, we get $\innerprthree{\Phi}{\Psi}\geqslant 0\,\forall_{\Psi\in\dual{\mappingcone}}$, which means that $\Phi\in\ddual{\mappingcone}$. But $\ddual{\mappingcone}=\mappingcone$ because $\mappingcone$ is a closed convex cone and Property \ref{propertyfour} holds. Hence $\Phi\in\mappingcone$. The proof of $1\Rightarrow 2$ strongly builds on the assumption that $\mappingcone$ has the mapping cone symmetry. By Proposition \ref{dualisamappingcone}, we know that $\dual{\mappingcone}$ is an mcs-cone as well. Therefore $\Psi\circ\Ad_V\in\dual{\mappingcone}$ for an arbitrary $\Psi\in\dual{\mappingcone}$ and $V\in\bk$. We have $\innerprthree{\Psi\circ\Ad_V}{\Phi}\geqslant 0\,\forall_{V\in\bk}\forall_{\Psi\in\dual{\mappingcone}}$. By Proposition \ref{propinnerpr}, point one, we get $\innerprthree{\Psi\circ{\Ad}_V}{\Phi}=\innerprthree{{\Ad}_V}{\conj{\Psi}\circ\Phi}$. Using Property \ref{propertythree} and Lemma \ref{lemmaCofAdV} with $\hilbertspacetwo=\hilbertspaceone$, the last term can be rewritten as
\begin{equation}\label{fourequalitiesforinnerpr}
 \innerprthree{{\Ad}_V}{\conj{\Psi}\circ\Phi}=\innerprtwo{\Choimatr{\Ad_V}}{\Choimatr{\conj{\Psi}\circ\Phi}}=\innerprtwo{\diad{v}}{\Choimatr{\conj{\Psi}\circ\Phi}}=\innerpr{\upsilon}{\Choimatr{\conj{\Psi}\circ\Phi}\left(\upsilon\right)},
\end{equation}
where $\upsilon=\sum_{i,j=1}^m V_{ij}f_j\otimes f_i$ for $V:\hilbertspaceone\ni a\mapsto\sum_{i,j=1}^mV_{ij}\innerpr{a}{f_j}f_i\in\hilbertspaceone$. The vector $\upsilon\in\hilbertspaceone\otimes\hilbertspaceone$ can be arbitrary, since we do not assume anything about the operator $V$. Consequently, the condition $\innerprthree{\Psi\circ\Ad_V}{\Phi}\geqslant 0\,\forall_{V\in\bk}\forall_{\Psi\in\dual{\mappingcone}}$ is equivalent to
\begin{equation}\label{posChoimatr}
 \innerpr{\upsilon}{\Choimatr{\conj{\Psi}\circ\Phi}\left(\upsilon\right)}\geqslant 0\,\forall_{\upsilon\in\hilbertspaceone\otimes\hilbertspaceone}\,\forall_{\Psi\in\dual{\mappingcone}},
\end{equation}
which means that $\Choimatr{\conj{\Psi}\circ\Phi}\in\boundedplus{\hilbertspaceone\otimes\hilbertspaceone}$ for all $\Psi\in\dual{\mappingcone}$. By the Choi theorem on completely positive maps \cite{ref.Choi75}, $\conj{\Psi}\circ\Phi\in\CPmapsb{\bk}$ for all $\Psi\in\dual{\mappingcone}$. Thus we have finished proving that $1\Leftrightarrow 2$. The proof of the equivalence $1\Leftrightarrow 3$ only needs a minor modification of the above argument. Instead of using point one of Proposition \ref{propinnerpr}, point two of the same proposition has to be used. Other details are practically the same as above and we shall not give them explicitly.
\end{proof}
\end{theorem}
In case of $\hilbertspacetwo=\hilbertspaceone$ and a $\conjsymb$-invariant mcs-cone $\mappingcone\in\Pmapsb{\bk}$, Theorem~\ref{MainTheorem} can be further simplified. 
\begin{theorem}\label{MainTheoremconjugatesymmetric}
Let $\mappingcone\subset\Pmapsb{\bk}$ be a $\conjsymb$-invariant mcs-cone. Then the following conditions are equivalent,
\begin{enumerate}
 \item $\Phi\in\mappingcone$,
\item $\Psi\circ\Phi\in\CPmapsb{\bk}$ for all $\Psi\in\dual{\mappingcone}$,
\item $\Phi\circ\Psi\in\CPmapsb{\bk}$ for all $\Psi\in\dual{\mappingcone}$.
\end{enumerate}
\begin{proof}
 Obvious from Theorem \ref{MainTheorem}.
\end{proof}
\end{theorem}
This result was earlier known for $\kPmapsb{k}{\bk}$ and $\kSPmapsb{k}{\bk}$ \cite{ref.SSZ09}, and inexplicitly for all symmetric (and convex) mapping cones \cite{ref.St09mappingcones}. As it was pointed to the author by Erling St{\o}rmer, in the case of $k$-positive maps, not necessarily from $\bk$  into itself, an even stronger characterization of the type of Theorems \ref{MainTheorem} and \ref{MainTheoremconjugatesymmetric} is valid. First, we have the simple
\begin{theorem}\label{MainTheoremkpositivemaps1}The following conditons are equivalent
 \begin{enumerate}
  \item $\Phi\in\kPmapsbb{k}{\bk}{\bh}$,
\item $\Ad_{\conj{V}}\circ\,\Phi\in\CPmapsb{\bk}$ for all $V\in\bbounded{\hilbertspaceone}{\hilbertspacetwo}$ s.t. $\rk V\leqslant k$,
\item $\Phi\circ\Ad_{\conj{V}}\in\CPmapsb{\bh}$ for all $V\in\bbounded{\hilbertspaceone}{\hilbertspacetwo}$ s.t. $\rk V\leqslant k$.
 \end{enumerate}
\begin{proof}
 Obvious from Theorem \ref{MainTheorem}. The duality relation 
\begin{equation}
 \dual{\kPmapsbb{k}{\bk}{\bh}}=\kSPmapsbb{k}{\bk}{\bh}=\convhull\left\{{\Ad}_V|V\in\bbounded{\hilbertspaceone}{\hilbertspacetwo},\rk V\leqslant k\right\}
\end{equation}
holds (cf. \cite{ref.SSZ09}) and we can substitute $\Psi$ in Theorem \ref{MainTheorem} with $\Ad_V$, $\rk V\leqslant k$. We also use the elementary fact that $\conj{\Ad_V}=\Ad_{\conj{V}}$.
\end{proof}
\end{theorem}
The next result on $k$-positive maps seems to be less obvious.
\begin{theorem}\label{MainTheoremkpositivemaps2}Denote with $\Pi_k\left(\hilbertspaceone\right)$ and $\Pi_k\left(\hilbertspacetwo\right)$ the sets of $k$-dimensional projections in $\hilbertspaceone$ and $\hilbertspacetwo$, resp.
 The following conditons are equivalent
 \begin{enumerate}
  \item $\Phi\in\kPmapsbb{k}{\bk}{\bh}$,
\item $\Ad_E\circ\,\Phi\in\CPmapsbb{\bk}{\bh}$ for all $E\in\Pi_k\left(\hilbertspacetwo\right)$,
\item $\Phi\circ\Ad_F\in\CPmapsbb{\bk}{\bh}$ for all $F\in\Pi_k\left(\hilbertspaceone\right)$,
\item $\Ad_E\circ\,\Phi\circ\Ad_F\in\CPmapsbb{\bk}{\bh}$ for all  $E\in\Pi_k\left(\hilbertspacetwo\right)$, $F\in\Pi_k\left(\hilbertspaceone\right)$.
 \end{enumerate}
\begin{proof}
 We shall prove the equivalence $1\Leftrightarrow 4$. The other ones follow analogously. Since $\dual{\CPmaps}=\CPmaps$ and any $\CPmaps$ map can be written as $\sum_i\Ad_{V_i}$ with $V_i$ arbitrary, the condition $\Ad_E\circ\,\Phi\circ\Ad_F\in\CPmapsbb{\bk}{\bh}$ is equivalent to
\begin{equation}\label{equivalentoffour}
 \innerprthree{{\Ad}_E\circ\,\Phi\circ{\Ad}_F}{{\Ad}_V}\geqslant 0\,\forall_{E\in\Pi_k\left(\hilbertspacetwo\right),F\in\Pi_k\left(\hilbertspaceone\right)}\forall_{V\in\bbounded{\hilbertspaceone}{\hilbertspacetwo}}.
\end{equation}
By Proposition \ref{propinnerpr}, point three, equation \eqref{equivalentoffour} can be rewritten as
\begin{equation}\label{equivalentoffour2}
 \innerprthree{\Phi}{{\Ad}_{EVF}}\geqslant 0\,\forall_{E\in\Pi_k\left(\hilbertspacetwo\right),F\in\Pi_k\left(\hilbertspaceone\right)}\forall_{V\in\bbounded{\hilbertspaceone}{\hilbertspacetwo}},
\end{equation}
where we used the fact that $\Ad_E\circ\Ad_V\circ\Ad_F=\Ad_{EVF}$ and the self-adjointness of $E$ and $F$. Note that $U=EVF$ is an element of $\bbounded{\hilbertspaceone}{\hilbertspacetwo}$ of rank $\leqslant k$. Conversely, every map in $U\in\bbounded{\hilbertspaceone}{\hilbertspacetwo}$ of rank $\leqslant k$ can be written in the form $EVF$ for some $V\in\bbounded{\hilbertspaceone}{\hilbertspacetwo}$, $E\in\Pi_k\left(\hilbertspacetwo\right)$ and $F\in\Pi_k\left(\hilbertspaceone\right)$. It is sufficient to take $V=U$ and $E$, $F$ as the range and rank projections for $U$, resp. Therefore the condition \eqref{equivalentoffour2} is equivalent to $\innerprthree{\Phi}{\Ad_U}\geqslant 0$ for all $U\in\bbounded{\hilbertspaceone}{\hilbertspacetwo}$ s.t. $\rk U\leqslant 0$. But this is the same as $\innerprthree{\Phi}{\Psi}\geqslant 0$ for all $\Psi\in\kSPmapsbb{k}{\bk}{\bh}$, or $\Phi\in\dual{\kSPmapsbb{k}{\bk}{\bh}}=\kPmapsbb{k}{\bk}{\bh}$. Thus $1\Leftrightarrow 4$.
\end{proof}
\end{theorem}

\begin{example}\label{exampleone}\textnormal{
A very instructive application of Theorem \ref{MainTheoremkpositivemaps2}, due to St{\o}rmer, has recently been given in \cite{ref.SS10}. It concerns maps of the form $\Phi_{\lambda}:\rho\mapsto\Tr\rho\cdot\Id-\lambda\Ad_V\left(\rho\right)$, or $\Phi_{\lambda}=\Tr-\lambda\Ad_V$ for short, where $V\in\bbounded{\hilbertspaceone}{\hilbertspacetwo}$ and $\lambda>0$. Consider first the question of complete positivity of such maps. It is not difficult to check that $\Choimatr{\Tr}=\One$. Thus by Lemma \ref{lemmaCofAdV}, $\Choimatr{\Phi_{\lambda}}=\One-\lambda\diad{\upsilon}$, where $\upsilon=\sum_{i=1}^n\sum_{j=1}^mV_{ij}f_j\otimes e_i$, $V:\hilbertspaceone\ni a\mapsto\sum_{i=1}^n\sum_{j=1}^mV_{ij}\innerpr{a}{f_j}e_i\in\hilbertspacetwo$ and $\One$ denotes the identity operator in $\hilbertspaceone\otimes\hilbertspacetwo$. According to the Choi theorem on completely positive maps \cite{ref.Choi75}, complete positivity of $\Phi_{\lambda}$ is equivalent to positivity of $\Choimatr{\Phi_{\lambda}}$. For an arbitrary vector $w\in\hilbertspaceone\otimes\hilbertspacetwo$, the product $\innerpr{w}{\Choimatr{\Phi_{\lambda}}\left(w\right)}$ equals $\innerpr{w}{w}-\lambda\left|\innerpr{w}{\upsilon}\right|^2$. Minimizing over vectors $w$ of unit norm, we get $1-\lambda\left|\innerpr{\upsilon}{\upsilon}\right|$. But $\left|\innerpr{\upsilon}{\upsilon}\right|=\sum_{i=1}^n\sum_{j=1}^mV_{ij}\overline{V_{ij}}$, which is the same as $\Tr\left(V\conj{V}\right)$. Thus $\Choimatr{\Phi_{\lambda}}$ is positive, or equivalently, $\Phi_{\lambda}$ completely positive if and only if $\lambda\Tr\left(V\conj{V}\right)\leqslant 1$. Now we turn to the question about $k$-positivity of $\Phi_{\lambda}$. According to Theorem \ref{MainTheoremkpositivemaps2}, condition 2., we need to check complete positivity of the map $\Ad_E\circ\,\Phi_{\lambda}$ for all projections $E$ of rank $k$. Interestingly, $\Ad_E\circ\,\Phi_{\lambda}=E\,\Tr-\lambda\Ad_{EV}=E\left(\Tr-\lambda\Ad_{EV}\right)E$, which is the same as $\Phi_{\lambda}$ save the $E$ at both ends and $EV$ in place of $V$. Actually, if we consider $\Ad_E\circ\,\Phi_{\lambda}$ as a map from $\bk$ into $\bounded{E\hilbertspacetwo}$, it is of the same form as $\Phi_{\lambda}$, with $EV$ instead of $V$. Hence, using the result obtained above, $\Ad_E\circ\,\Phi_{\lambda}$ is a $\CPmaps$ map if and only if $\lambda\Tr\left(EV\conj{\left(EV\right)}\right)=\lambda\Tr\left(EV\conj{V}\right)\leqslant 1$, where we skip a few technical details of the argument (cf. \cite{ref.SS10}). According to Theorem \ref{MainTheoremkpositivemaps2}, we should maximize the expression on the left over all possible choices of $E$. The maximum turns out to be equal to $\lambda$ times the square of the $k$-th Ky Fan norm of $V$ (cf. e.g. \cite{ref.HornJohnson}). Thus we rederive a result by Chru\'sci\'nski and Kossakowski concerning $k$-positivity of maps of the form $\Phi_{\lambda}$ \cite{ref.ChK09}.} 
\end{example}
  
\section{Conclusion}
In the finite-dimensional setting,
we presented a number of results concerning convex cones with a mapping cone symmetry, or ``mcs-cones''. Our focus was on convex cone duality. The use of a slightly modified class of cones, but much in the spirit of the original definition of a mapping cone by St\o rmer \cite{ref.St86}, allowed us to make very general statements. In Proposition~\ref{dualisamappingcone}, we showed that the dual of an mcs-cone is an mcs-cone.
Our main result, which is the characterization included in Theorem \ref{MainTheorem}, can be very loosely described as saying that the surface of mcs-cones has an additional structure, which makes them more ``smooth''. There is a stronger relation between a pair of dual cones \emph{with a mapping cone symmetry} than the relation between a mere pair of dual convex cones. There also exist stronger versions of the main characterization theorem, valid for $\conjsymb$-invariant cones (Theorem \ref{MainTheoremconjugatesymmetric}) and specifically for $k$-positive maps (Theorems \ref{MainTheoremkpositivemaps1} and \ref{MainTheoremkpositivemaps2}). Example~\ref{exampleone} shows a practical application. 

\begin{figure}\begin{center}
               \includegraphics[scale=0.33]{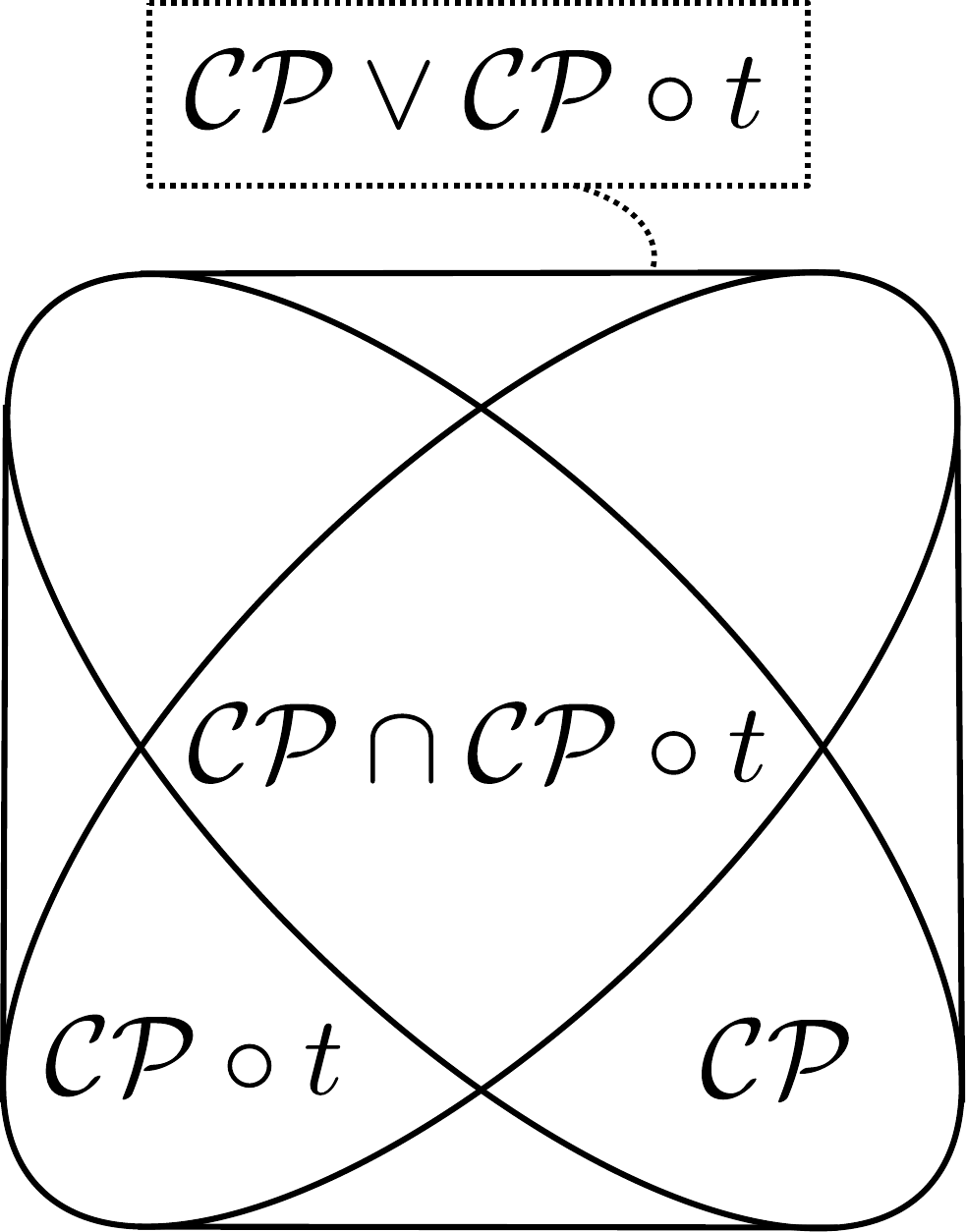}
              \end{center} 
\caption{\label{figureCPCPt}A family of four mapping cones, constructed out of the cone of completely positive maps. The set $\CPmaps\lor\CPmaps\circ t$ of decomposable maps is included.}
\end{figure}

It is natural to ask how large the class of mcs-cones is. We know that $k$-positive and $k$-superpositive maps (for $k=1,2,\ldots$) provide examples of such cones, including completely positive maps. Another class includes the same cones multiplied by the transposition map, i.e. cones $\mappingcone\circ\transpose:=\left\{\Phi\circ\transpose\,\vline\,\Phi\in\mappingcone\right\}$, where $\mappingcone$ stands for any of $\Pmaps$, $\kPmaps{k}$, $\CPmaps$, $\kSPmaps{k}$, $\SPmaps$. For example, elements of $\CPmaps\circ\transpose$ are called {\it completely co-positive}. We can also provide a variety of new examples by taking the intersection and the convex sum of any two known mapping cones $\mappingcone_1$ and $\mappingcone_2$, one not included in another. Figure \ref{figureCPCPt} shows roughly how it works for $\mappingcone_1=\CPmaps$ and $\mappingcone_2=\CPmaps\circ\transpose$. 
The key remaining question is 
\begin{question}\label{questionone}
 Are there any ``untypical'' mcs-cones?
\end{question}
By a typical mcs-cone we mean a cone which is obtained from $\Pmaps$, $\kPmaps{k}$, $\CPmaps$, $\kSPmaps{k}$ or $\SPmaps$ using the mentioned operations $\mappingcone\rightarrow\mappingcone\circ\transpose$, $\left(\mappingcone_1,\mappingcone_2\right)\rightarrow\mappingcone_1\cap\mappingcone_2$ and $\left(\mappingcone_1,\mappingcone_2\right)\rightarrow\mappingcone_1\lor\mappingcone_2$. So far, no answer to the question is known. Note that all typical mcs-cones are symmetric.

At this place, it seems desirable to shortly mention a connection between mcs-cones and matrix ordered  $\ast$-vector spaces \cite{ref.ChEff77}, or operator systems. For any mcs-cone $\mappingcone\in\bounded{\bk}$, one has a matrix ordering of $\bk$, given by the cones $C^{\mappingcone}_n:=\left\{\left(X_{ij}\right)\in\matrices{n}{\bk}\vline\left(\Phi\left(X_{ij}\right)\right)\in\matrices{n}{\bk}^+\forall_{\Phi\in\mappingcone}\right\}$ for $n=1,2,\ldots$ . A similar definition could also be used in the case when $\mappingcone$ does not have a mapping cone symmetry, however mcs-cones seem to be a preferred choice, as they provide a definition of positivity of $\left(X_{ij}\right)\in\matrices{n}{\bk}$ which is invariant under any conjugation map $\left(X_{ij}\right)\mapsto\left(\Ad_V\left(X_{ij}\right)\right)=\left(VX_{ij}\conj{V}\right)$. This is particularly natural if we think about the case $n=1$. The mapping cone symmetry also takes care of a redundancy in the condition $\left(\Phi\left(X_{ij}\right)\right)\in\matrices{n}{\bk}^+\forall_{\Phi\in\mappingcone}$ allowing a conjugation of $\Phi$, $\left(\Ad_U\circ\Phi\left(X_{ij}\right)\right)\in\matrices{n}{\bk}^+\forall_{\Phi\in\mappingcone}\forall_{U\in\bk}$, to have no effect on $C^{\mappingcone}_n$. Thus a left-invariance of $\mappingcone$ under conjugation can always be assumed when a matrix ordering given by the cones $C^{\mappingcone}_n$ for some $\mappingcone\subset\bounded{\bk}$ is considered. In many respects, orderings of this type provide a generalization of the $\textnormal{OMIN}^k$ and $\textnormal{OMAX}^k$ structures recently investigated in \cite{ref.JKPP10}. In that broader context, Question \ref{questionone} alludes to the subject of possible operator system structures for $\bk$.

Throughout the paper, we used the definition \eqref{HSProdtwo} of an inner product in the space of linear maps from an algebra $\bk$ into $\bh$, where $\hilbertspaceone$, $\hilbertspacetwo$ stand for two Hilbert spaces and we assumed $\hilbertspaceone$ and $\hilbertspacetwo$ to be finite-dimensional. We also exploited the property that the Jamio{\l}kowski-Choi isomorphism (cf. eq. \eqref{Jamisodef}) is an isometry. Because of the finite-dimensionality assumption, the definition \eqref{HSProdtwo} was certainly correct and quite natural. The same for the Jamio{\l}kowski-Choi isomorphism. If the dimension of $\hilbertspaceone$ or $\hilbertspacetwo$ was not assumed to be finite, there will be serious problems with both definitions. Nevertheless, it seems that similar methods may work at least for some cones in the infinite-dimensional case, e.g. assuming that their elements are trace class in a proper sense.

In the end, let us mention that $k$-positive maps are of special interest in the theory of quantum information and computation, a very 
active branch of contemporary physics and information science (cf. e.g. \cite{ref.Kribs05}). The case $k=1$ corresponds to entanglement witnesses \cite{ref.HHH96} and $k=2$ to undistillable quantum states, with the fundamental question about the existence of NPT bound entanglement \cite{ref.DiVicenzo00}.

\section{Acknowledgement}
The author would like to thank Erling St\o rmer, Marcin Marciniak and Karol \.Zyczkowski for collaboration and for comments on the manuscript. The project was operated within the Foundation
for Polish Science International Ph.D. Projects Programme co-financed
by the European Regional Development Fund covering, under the agreement
no. MPD/2009/6, the Jagiellonian University International Ph.D. Studies in
Physics of Complex Systems. Some of the results presented in the paper came to existence during a visit of the author to the University of Oslo, where he enjoyed hospitality at the Mathematics Institute and was financially supported by Scholarschip and Training Fund, operated by Foundation for the Development of the Education System.

\bibliographystyle{elsarticle-num-names}
\bibliography{Mapping_cones_LAA_corrected}

\begin{thebibliography}{19}
\providecommand{\natexlab}[1]{#1}
\providecommand{\url}[1]{\texttt{#1}}
\providecommand{\urlprefix}{URL }
\expandafter\ifx\csname urlstyle\endcsname\relax
  \providecommand{\doi}[1]{doi:\discretionary{}{}{}#1}\else
  \providecommand{\doi}[1]{doi:\discretionary{}{}{}\begingroup
  \urlstyle{rm}\url{#1}\endgroup}\fi
\providecommand{\bibinfo}[2]{#2}

\bibitem[{St{\o}rmer(1986)}]{ref.St86}
\bibinfo{author}{E.~St{\o}rmer}, \bibinfo{title}{Extension of positive maps
  into $\bh$}, \bibinfo{journal}{J. Funct. Anal.} \bibinfo{volume}{66}
  (\bibinfo{year}{1986}) \bibinfo{pages}{235--254}.

\bibitem[{St{\o}rmer(2008)}]{ref.St08}
\bibinfo{author}{E.~St{\o}rmer}, \bibinfo{title}{Separable states and positive
  maps}, \bibinfo{journal}{J. Funct. Anal.} \bibinfo{volume}{254}
  (\bibinfo{year}{2008}) \bibinfo{pages}{2303--2312}.

\bibitem[{St{\o}rmer(2009{\natexlab{a}})}]{ref.St09dual}
\bibinfo{author}{E.~St{\o}rmer}, \bibinfo{title}{Duality of cones of positive
  maps}, \bibinfo{journal}{M{\"u}nster J. Math.} \bibinfo{volume}{2}
  (\bibinfo{year}{2009}{\natexlab{a}}) \bibinfo{pages}{299--310}.

\bibitem[{St{\o}rmer(2009{\natexlab{b}})}]{ref.St09sepposII}
\bibinfo{author}{E.~St{\o}rmer}, \bibinfo{title}{Separable states and positive
  maps {II}}, \bibinfo{journal}{Math. Scand.} \bibinfo{volume}{105}
  (\bibinfo{year}{2009}{\natexlab{b}}) \bibinfo{pages}{188--198}.

\bibitem[{Skowronek et~al.(2009)Skowronek, St{\o}rmer, and
  {\.Z}yczkowski}]{ref.SSZ09}
\bibinfo{author}{{\L}.~Skowronek}, \bibinfo{author}{E.~St{\o}rmer},
  \bibinfo{author}{K.~{\.Z}yczkowski}, \bibinfo{title}{Cones of positive maps
  and their duality relations}, \bibinfo{journal}{J. Math. Phys.}
  \bibinfo{volume}{50} (\bibinfo{year}{2009}) \bibinfo{pages}{062106}.

\bibitem[{Choi(1975)}]{ref.Choi75}
\bibinfo{author}{M.-D. Choi}, \bibinfo{title}{Completely positive linear maps
  on complex matrices}, \bibinfo{journal}{Lin. Alg. Appl.} \bibinfo{volume}{10}
  (\bibinfo{year}{1975}) \bibinfo{pages}{285--290}.

\bibitem[{Jamio{\l}kowski(1972)}]{ref.J72}
\bibinfo{author}{A.~Jamio{\l}kowski}, \bibinfo{title}{Linear transformations
  which preserve trace and positive semidefinitness of operators},
  \bibinfo{journal}{Rep. Math. Phys.} \bibinfo{volume}{3}
  (\bibinfo{year}{1972}) \bibinfo{pages}{275--278}.

\bibitem[{Ando(2004)}]{ref.Ando04}
\bibinfo{author}{T.~Ando}, \bibinfo{title}{Cones and norms in the tensor
  product of matrix spaces}, \bibinfo{journal}{Lin. Alg. Appl.}
  \bibinfo{volume}{379} (\bibinfo{year}{2004}) \bibinfo{pages}{3--41}.

\bibitem[{De~Pillis(1967)}]{ref.Pillis}
\bibinfo{author}{J.~De~Pillis}, \bibinfo{title}{Linear transformations which
  preserve {H}ermitian and positive semidefinite operators},
  \bibinfo{journal}{Pacific J. Math.} \bibinfo{volume}{23}
  (\bibinfo{year}{1967}) \bibinfo{pages}{129--137}.

\bibitem[{Rockafellar(1997)}]{ref.Rockafellar}
\bibinfo{author}{R.~Rockafellar}, \bibinfo{title}{Convex {A}nalysis},
  \bibinfo{publisher}{Princeton University Press}, \bibinfo{year}{1997}.

\bibitem[{St{\o}rmer(2010)}]{ref.St09mappingcones}
\bibinfo{author}{E.~St{\o}rmer}, \bibinfo{title}{Mapping cones of positive
  maps}, \bibinfo{journal}{Math. Scand.} \bibinfo{note}{To appear, preprint
  arXiv:0906.0472}.

\bibitem[{Skowronek and St{\o}rmer(2010)}]{ref.SS10}
\bibinfo{author}{{\L}.~Skowronek}, \bibinfo{author}{E.~St{\o}rmer},
  \bibinfo{title}{Choi matrices, norms and entanglement associated with
  positive maps on matrix algebras}, \bibinfo{note}{preprint arXiv:1008.3126},
  \bibinfo{year}{2010}.

\bibitem[{Horn and Johnson(1991)}]{ref.HornJohnson}
\bibinfo{author}{R.~Horn}, \bibinfo{author}{C.~Johnson}, \bibinfo{title}{Topics
  in {M}atrix {A}nalysis}, \bibinfo{publisher}{Cambridge University Press},
  \bibinfo{year}{1991}.

\bibitem[{Chru{\'s}ci{\'n}ski and Kossakowski(2009)}]{ref.ChK09}
\bibinfo{author}{D.~Chru{\'s}ci{\'n}ski}, \bibinfo{author}{A.~Kossakowski},
  \bibinfo{title}{Spectral conditions for positive maps},
  \bibinfo{journal}{Commun. Math. Phys.} \bibinfo{volume}{290}
  (\bibinfo{year}{2009}) \bibinfo{pages}{1051--1064}.

\bibitem[{Choi and Effros(1977)}]{ref.ChEff77}
\bibinfo{author}{M.-D. Choi}, \bibinfo{author}{E.~G. Effros},
  \bibinfo{title}{Injectivity and operator spaces}, \bibinfo{journal}{J. Funct.
  Anal.} \bibinfo{volume}{24} (\bibinfo{year}{1977}) \bibinfo{pages}{156--209}.

\bibitem[{Johnston et~al.(2010)Johnston, Kribs, Paulsen, and
  Pereira}]{ref.JKPP10}
\bibinfo{author}{N.~Johnston}, \bibinfo{author}{D.~W. Kribs},
  \bibinfo{author}{V.~I. Paulsen}, \bibinfo{author}{R.~Pereira},
  \bibinfo{title}{Minimal and Maximal Operator Spaces and Operator Systems in
  Entanglement Theory}, \bibinfo{journal}{J. Funct. Anal.} \bibinfo{note}{To
  appear, preprint arXiv:1010.1432}.

\bibitem[{Kribs(2005)}]{ref.Kribs05}
\bibinfo{author}{D.~Kribs}, \bibinfo{title}{A quantum computing primer for
  operator theorist}, \bibinfo{journal}{Lin. Alg. Appl.} \bibinfo{volume}{400}
  (\bibinfo{year}{2005}) \bibinfo{pages}{147--167}.

\bibitem[{Horodecki et~al.(1996)Horodecki, Horodecki, and
  Horodecki}]{ref.HHH96}
\bibinfo{author}{M.~Horodecki}, \bibinfo{author}{P.~Horodecki},
  \bibinfo{author}{R.~Horodecki}, \bibinfo{title}{Separability of mixed states:
  necessary and sufficient conditions}, \bibinfo{journal}{Phys. Lett. A}
  \bibinfo{volume}{223} (\bibinfo{year}{1996}) \bibinfo{pages}{1--8}.

\bibitem[{DiVincenzo et~al.(2000)DiVincenzo, Shor, Smolin, Terhal, and
  Thapliyal}]{ref.DiVicenzo00}
\bibinfo{author}{D.~DiVincenzo}, \bibinfo{author}{P.~Shor},
  \bibinfo{author}{J.~Smolin}, \bibinfo{author}{B.~Terhal},
  \bibinfo{author}{A.~Thapliyal}, \bibinfo{title}{Evidence for bound entangled
  states with negative partial transpose}, \bibinfo{journal}{Phys. Rev. A}
  \bibinfo{volume}{61} (\bibinfo{year}{2000}) \bibinfo{pages}{062312}.

\end{thebibliography}

\end{document}